\documentclass[11pt]{article}\usepackage{graphicx}\usepackage[T1]{fontenc}
\usepackage[utf8]{inputenc}\usepackage{xcolor}
\usepackage[dvips]{epsfig}
\usepackage{amssymb}
\newtheorem{theorem}{Theorem} \newtheorem{lemma}{Lemma} 
\newtheorem{remark}{Remark}\newtheorem{proposition}{Proposition}
\newtheorem{corollary}{Corollary} \newcommand{\La}{\Lambda}
\newcommand{\R}{{\mathbb R}}  \newcommand{\Z}{{\mathbb Z}} \newcommand{\N}{{\mathbb N}}
 \newcommand{\C}{{\mathbb C}}  

\begin{document}\title{Fourier quasicrystals  with unit masses}

\author{Alexander Olevskii  and  Alexander Ulanovskii}

\date{}\maketitle

    \begin{abstract}Every  set $\La\subset\R$ such that the sum of $\delta$-measures  sitting at the points of $\La$
                        is a  Fourier quasicrystal, is  the zero set of an exponential polynomial with imaginary frequencies.\end{abstract}


\section{Introduction}

By a {\it crystalline measure} one means an atomic measure $\mu$ which is a tempered distribution and whose distributional Fourier transform is also an atomic measure, $$
               \mu=\sum_{\lambda\in\La}c_\lambda\delta_\lambda,\quad \hat\mu=\sum_{s\in S}a_s\delta_s,\quad c_\lambda,a_s\in \C,
               $$where   the support $\La$ and  the spectrum $S$  of  $\mu$  are locally finite sets. If in addition the measures  $|\mu|$ and $|\hat\mu|$ are also tempered, then $\mu$ is called a {\it Fourier quasicrystal} (FQ). A
                classical example of an FQ  is a Dirac comb
                $$
                \mu=\sum_{n\in\Z}\delta_n,
                $$which satisfies $\hat\mu=\mu$. Here $\delta_x$  denotes  the Dirac measure at point $x$.

                A few constructions of {\it aperiodic}  FQs are known, see \cite{lo}, \cite{k16}, \cite{m16}.

                 Recently P. Kurasov and P. Sarnak \cite{ks} found  examples of FQs  with unit masses,
\begin{equation}\label{mu}
\mu=\sum_{\lambda\in\La}\delta_\lambda,
\end{equation}
         where $\La$ is an aperiodic {\it uniformly discrete set} in $\R$.

                  Alternative approaches to construction of such measures were suggested by Y.~Meyer \cite{m20} and later in \cite{ou}.   In the present note we
                show that the  construction  in the last paper characterises
                all FQs of the form (\ref{mu}):

\begin{theorem}   Let $\mu$  be  an FQ of the form (\ref{mu}). 
                          Then there is an exponential  polynomial $$  p(x)=\sum_{k=1}^nb_ke^{i\gamma_k x},\quad n\in\N, b_k\in\C, \gamma_k\in\R, k=1,...,n,$$  with real simple zeros
                           such that $\La$ is the zero set of $p$.\end{theorem}

In the opposite direction, the construction in \cite{ou} shows that for every exponential polynomial $p$ with real simple zeros there is an FQ $\mu$ with unit masses  whose support $\La$ is the zero set of $p$.

\section{Auxiliary Results}

We will use the standard form of Fourier transform
$$
\hat h(t)=\int_\R e^{-2\pi itx}h(x)\,dx,\quad h\in L^1(\R).
$$

Let us start with a result which may have intrinsic interest:

\begin{proposition} Let $\mu$ be a positive measure which is                a tempered distribution, such that its distributional Fourier transform $\hat \mu$ is  a measure satisfying \begin{equation}\label{2}
                                 \int_{(-R,R)}d|\hat\mu|  = O(R^m),\ R\to\infty, \quad   \mbox{for some} \ m>0,
\end{equation}which means that $|\hat\mu|$ is a tempered distribution.
Then there exists $C$ such that
\begin{equation}\label{1} \mu(a,b)\leq C(1+b-a),\quad  -\infty<a<b<\infty.\end{equation}\end{proposition}

\noindent{\bf Proof}.
It suffices to prove (\ref{1}) for every interval $(a,b)$ satisfying $b-a\geq 2.$

Fix any  non-negative Schwartz function $g(x)$ supported by $[-1/2,1/2]$ and such that
$$\int_\R g(x)\,dx=1.
$$Set
$$
f(x):=(g\ast 1_{(a-1/2,b+1/2)})(x)\in S(\R).
$$Clearly,
$$|\hat f(t)|=  |\hat g(t) \widehat{ 1}_{(a-1/2,b+1/2)}(t)|\leq (1+b-a)|\hat g(t)|.
$$
 Using this inequality and (\ref{2}), we get
$$
\int_\R f(x)\mu(dx)= \int_\R\hat f(t)\hat\mu(dt)\leq (1+b-a)\int_{\R} |\hat g(t)| |\hat\mu|(dt)=C(1+b-a).
$$

On the other hand, clearly,
$$
f(x)=g(x)\ast 1_{(a-1/2,b+1/2)}(x) =1, \quad x\in (a,b).
$$Hence,
$$
\int_\R f(t)\mu(dt)\geq \mu(a,b),
$$which proves the proposition.

Recall that a set $\La\subset\R$ is called uniformly discrete, if $$ \sup_{\lambda',\lambda\in\La,\lambda\ne\lambda'}|\lambda-\lambda'|>0.$$A set  $\La$ is called relatively uniformly discrete if it is a union of finite number of uniformly discrete sets. $\Box$

 Proposition 1 implies
\begin{corollary}Let $\mu$ be a  measure of the form (\ref{mu}) whose distributional Fourier transform is a measure satisfying (\ref{2}).  Then its support $\La$ is a relatively uniformly discrete set.
\end{corollary}

In what follows we assume that $\mu$ is an FQ of the form (\ref{mu}). Recall that  its Fourier transform $\hat\mu$
admits a representation
\begin{equation}\label{hmu}
\hat\mu=\sum_{s\in S}a_s\delta_s,
\end{equation}where $S$ is a locally finite set and $|\hat\mu|$ satisfies (\ref{2}).

We  assume that $0\not\in\La$ (otherwise, we consider the measure $\mu(d(x-a))$, for an appropriate $a$).

\begin{lemma}  The formula\begin{equation}\label{f}
\sum_{\lambda\in\La}\left(\frac{1}{w-\lambda}+\frac{1}{\lambda}\right)=\alpha-2\pi i\sum_{s\in S,s<0}a_s e^{-2\pi iws},\quad w\in\C, \mbox{\rm Im}\,w>0,\end{equation}holds with some $\alpha\in\C$.\end{lemma}

\noindent{\bf Proof}.  For every $w\in\C$, Im$\,w\geq1,$ set
$$
e_w(t)=\left\{\begin{array}{ll} 2\pi e^{-2\pi iwt}&t<0\\ 0&t\geq0 \end{array}\right.
$$Then the inverse Fourier transform of $e_w$ is given by $$\check{e}_w(x)=\frac{i}{w-x},\quad x\in\R.$$

Denote by $s_1\in S, s_1<0,$  the greatest negative element of $S$, were $S$ is the spectrum of $\hat\mu$ in (\ref{hmu}).
Choose any non-negative function $H\in S(\R)$ which vanishes outside of  $(-1,0)$  and such that
$$
\int_\R H(t)\,dt=1.
$$Also, denote by $h$ the inverse Fourier transform of $H$.

 Fix $\epsilon, 0<\epsilon<-s_1$, and set
$$
f_\epsilon(x):=\frac{ih(\epsilon x)}{w-x}.
$$  Then its Fourier transform is given by
$$
\hat f_\epsilon(t):=\frac{1}{\epsilon}e_w(t)\ast H(t/\epsilon).
$$

The following properties of $\hat f_\epsilon$ are evident: $\hat f_\epsilon\in S(\R)$, $f_\epsilon(t)=0,t\geq0$ and
$$
\hat f_\epsilon(t)=\frac{2\pi}{\epsilon}\int_{-\epsilon}^0 e^{-2\pi i w(t-u)}H(u/\epsilon)\,du=2\pi h(\epsilon w)e^{-2\pi iwt},\quad t\leq s_1<- \epsilon.
$$

Hence, by (\ref{hmu}),
$$
\sum_{\lambda\in\La}f_\epsilon(\lambda)=\sum_{\lambda\in\La}\frac{ih(\epsilon \lambda)}{w-\lambda}=2\pi h(\epsilon w)\sum_{s\in S,s<0}a_se^{-2\pi iws}.
$$The series above converges absolutely due to (\ref{2}).
Letting $w=i$, we get
$$
\sum_{\lambda\in\La}\frac{ih(\epsilon \lambda)}{i-\lambda}= h(\epsilon i)\alpha_0,
$$where
$$
\alpha_0:=2\pi \sum_{s\in S,s<0}a_se^{2\pi s}.
$$
We conclude that
$$
\sum_{\lambda\in\La}h(\epsilon \lambda)\left(\frac{1}{w-\lambda}-\frac{1}{i-\lambda}\right)=ih(\epsilon i)\alpha_0- 2\pi ih(\epsilon w)\sum_{s\in S,s<0}a_se^{-2\pi iws}.
$$

Now,  by Corollary 1 the series converges
$$
\sum_{\lambda\in\La}\left|\frac{1}{w-\lambda}-\frac{1}{i-\lambda}\right|<\infty,\quad \mbox{\rm Im}\,w>0.
$$Hence,  in the formula above we may let $\epsilon\to0$ to get
$$
\sum_{\lambda\in\La}\left(\frac{1}{w-\lambda}-\frac{1}{i-\lambda}\right)=i\alpha_0-2\pi i\sum_{s\in S,s<0}a_se^{-2\pi iws}.
$$
This proves  (\ref{f}), where $\alpha$ is defined as
$$
\alpha:=i\alpha_0+\sum_{\lambda\in\La}\left(\frac{1}{\lambda}+\frac{1}{i-\lambda}\right).
$$The series above converges absolutely due to Corollary 1.
$\Box$

Similarly to the proof above, one can establish a variant of (\ref{f}) for the lower half-plane:
\begin{equation}\label{ff}\sum_{\lambda\in\La}\left(\frac{1}{w-\lambda}+\frac{1}{\lambda}\right)=\beta+2\pi i\sum_{s\in S,s>0}a_s e^{-2\pi iws},\quad w\in\C,\mbox{\rm Im}\,w<0,\end{equation}where $\beta\in\C$.

\section{Proof of Theorem 1}

 It follows from Corollary 1 that
$$
\sum_{\lambda\in\La,\lambda\ne0}\frac{1}{|\lambda|^{1+\epsilon}}<\infty,\quad \epsilon>0.
$$

Set
\begin{equation}\label{psi}
\psi(z):=\lim_{R\to\infty}\prod_{\lambda\in\La,|\lambda|<R}\left(1-\frac{z}{\lambda}\right)e^{z/\lambda}.
\end{equation}
The inequality above and a theorem of Borel (see \cite{le}, Lecture 4, Theorem~3) imply that $\psi$ is an entire function of order one, i.e. for every $\epsilon>0$ there is a constant $C(\epsilon)$ such that
$$
\max_{z:|z|\leq R}|\psi(z)|\leq C(\epsilon)e^{|z|^{1+\epsilon}}, \quad R>0.
$$

It is clear that we have
$$
\frac{\psi'(w)}{\psi(w)}=\sum_{\lambda\in\La}\left(\frac{1}{w-\lambda}+\frac{1}{\lambda}\right),\quad w\in\C.
$$
This and (\ref{f}) imply the representation
\begin{equation}\label{psi}
\psi(w)=Ce^{\alpha w}e^{\sum_{s\in S, s<0}(a_s/ s)e^{-2\pi isw}},\quad \mbox{{\rm Im}}\,w>0,\end{equation}with some $C\in\C.$

 Set
$$
p(w):=\psi(w)e^{-\alpha w}.
$$Since $\psi$ is of order one, the same is true for $p$.

Observe that  $p$ has simple real zero at the points of $\La.$
To prove the theorem, we show that $p$ is an exponential polynomial with imaginary frequencies.

By (\ref{psi}),
\begin{equation}\label{p}
p(w)=Ce^{\sum_{s\in S, s<0}(a_s/s)e^{-2\pi isw}},\quad \mbox{{\rm Im}}\,w>0.
\end{equation}
 This  shows  that $p$ is bounded from above in any  half-plane  Im$w\geq r, r>0$.
  Similarly, by  (\ref{ff}), the function $$
\psi(w)e^{-\beta w}=p(w)e^{(\alpha-\beta)w}
$$ is bounded from above in any half-plane  Im$w\leq -r, r>0.$
Hence,  $p$  admits an estimate
$$
|p(w)|\leq Ke^{|\alpha-\beta||w|},\quad \mbox{\rm Im}\,w\leq-1,\quad K>0.
$$

We now apply a variant of the Phragm\'{e}n--Lindel\"{o}f principle (see \cite{le}, Lecture 6, Theorem 2)  in the angles
$\{|$arg$\,w|\leq \pi/4\}$ and $\{|\pi-$arg$\,w|\leq \pi/4\}$  to prove that $p$ is an entire functions of exponential type, i.e. $$|p(w)|\leq Ce^{\sigma|w|},\quad w\in\C,$$with some $C,\sigma>0$.

 Recall that $p$ is bounded in the half-plane Im$\,w\geq 1$. Hence, another variant of the Phragm\'{e}n--Lindel\"{o}f principle (\cite{le}, Lecture 6, Theorem~3) implies that $p$  is bounded in every  strip $|$Im$\,w|\leq r, r>0,$ and in the lower half-plane it satisfies the estimate
\begin{equation}\label{pp}
|p(w)|\leq C e^{- \sigma y},\quad w=x+iy, -\infty<y\leq 0.
\end{equation}

 Now, develop each exponent in (\ref{p}) into a series
$$\exp((a_s/s)e^{-2\pi isw})               =\sum_{k=0}^\infty\left(\frac{a_s}{s}\right)^k e^{-2\pi i ksw},\quad \mbox{\rm Im}\,w>0,
$$
 to get the representation\begin{equation}\label{c}p(w)=\sum_{u\in U}b_ue^{-2\pi iuw},\quad \mbox{Im}\,w>0,\end{equation}
where $U$ is defined as follows
$$
U:= \{0\}\cup S_-\cup(S_-+S_-)\cup(S_-+S_-+S_-)\cup...\subset(-\infty,0),\quad S_-:=S\cap(-\infty,0).
$$
One may check that $U$ is a locally finite set and that the series in (\ref{c}) converges absolutely for every $w\in\C,$ Im$\,w>0.$

 To prove the theorem, we have to show that the series in the right hand-side of (\ref{c}) contains only a finite number of terms.

Set
$$
p_\epsilon(w):=p(w)\left(\frac{\sin\pi\epsilon w}{\pi\epsilon w}\right)^2.
$$
  Recall that  $p$ is bounded in the upper half-pane and satisfies  (\ref{pp}) in the lower one. It is easy to check that the function $p_\epsilon$ has the following properties: It is an entire function of exponential type, the inequalities  $$
  |p_\epsilon(w)|\leq Ce^{2\pi\epsilon y},\ \mbox{\rm Im}\, w\geq0,\quad |p_\epsilon(w)|\leq Ce^{(\sigma+2\pi\epsilon) |y|},\ \mbox{\rm Im}\, w\leq0, \ w=x+iy,
  $$hold in the upper and lower half-planes, and finally it belongs to $L^2\cap L^1$ on every line Im$\,w=const.$ Using the Paley--Wiener theorem, we see that  the inverse Fourier transform of $p_\epsilon$, $$\check{p}_\epsilon(t)=\int_\R e^{2\pi i tw}p_\epsilon(w)\,dx,\quad w=x+iy, $$ is a continuous function which vanishes
                outside the  interval $(-\sigma-\epsilon,\epsilon).$

 Let us check that in (\ref{c}) we have $b_u=0,u\in U,$ whenever $u<-\sigma$.
We choose $\epsilon>0$ so small that $u<-\sigma-\epsilon$ and the interval $(u-\epsilon,u+\epsilon)$ contains no other points of  $U$. Then, by (\ref{c}),
$$
0=\check{p}_\epsilon(u)=b_u\int_\R \left(\frac{\sin\pi\epsilon w}{\pi\epsilon w}\right)^2\,dx=b_u\int_\R \left(\frac{\sin\pi\epsilon x}{\pi\epsilon x}\right)^2\,dx,
$$which proves that $b_u=0.$
$\Box$

\begin{remark} 1. Some minor changes in the proof of Theorem 1 allow one to prove the following extension of
Theorem 1 to measures with integer masses,
\begin{equation}\label{mm}
\mu=\sum_{\lambda\in\La}c_\lambda\delta_\lambda,\quad c_\lambda\in\N:
\end{equation} If such a measure is an FQ, then there is an exponential polynomial $p(x)$ with imaginary frequencies and  (not necessarily simple) real zeros such that $\La$ is the zero set of $p$.

2. In the opposite direction, using the approach in \cite{ou}, one may check that to  every exponential polynomial $p(x)$ with imaginary frequencies and  real zeros  there corresponds  an FQ $\mu$ of the form (\ref{mm}), where $\La$ is the zero set of $p$.

\end{remark}

\end{document}